# UNIFORMLY ROOT-$N$ CONSISTENT DENSITY ESTIMATORS FOR WEAKLY DEPENDENT INVERTIBLE LINEAR PROCESSES

BY ANTON SCHICK[1] AND WOLFGANG WEFELMEYER

*Binghamton University and University of Cologne*

Convergence rates of kernel density estimators for stationary time series are well studied. For invertible linear processes, we construct a new density estimator that converges, in the supremum norm, at the better, parametric, rate $n^{-1/2}$. Our estimator is a convolution of two different residual-based kernel estimators. We obtain in particular convergence rates for such residual-based kernel estimators; these results are of independent interest.

**1. Introduction.** The usual estimators for the density of a stationary process are kernel estimators and their recursive versions. Rates of convergence and pointwise central limit theorems have been studied under various mixing conditions by Robinson [24], Chanda [8], Castellana and Leadbetter [7], Masry [19, 20, 21, 22], Tran [39, 40, 41], Roussas [27, 28, 29], Cai and Roussas [6], Ango Nze and Portier [2], Ango Nze and Doukhan [1], Ango Nze and Rios [3], Doukhan and Louhichi [11] and Dedecker and Merlevède [10], and for linear processes by Hall and Hart [14], Tran [42], Hallin and Tran [15], Coulon-Prieur and Doukhan [9], Honda [16], Lu [18], Wu and Mielniczuk [43], Bryk and Mielniczuk [5] and Schick and Wefelmeyer [36, 37]. Under appropriate conditions, the convergence rates of these kernel estimators are the same as for independent and identically distributed observations.

Linear processes are written as linear combinations of independent innovations and the stationary density can be represented as a convolution of other densities in many different ways. We use the simplest such representation and estimate the stationary density by plugging in residual-based estimators of the densities involved in the representation. We expect this

Received December 2003; revised June 2006.
[1]Supported in part by NSF Grant DMS-04-05791.
[AMS 2000 subject classifications.](...) 62G07, 62G20, 62M05, 62M10.

to lead to faster, parametric rates of convergence. This is already known in nonparametric models with i.i.d. observations. Frees [12] shows that his plug-in estimators for densities of certain functions $q(X_1, \ldots, X_m)$ are pointwise $n^{1/2}$-consistent. Saavedra and Cao [32] consider the special case $q(X_1, X_2) = X_1 + aX_2$. Schick and Wefelmeyer [34, 38] prove functional convergence for $q(X_1, \ldots, X_m) = u_1(X_1) + \cdots + u_m(X_m)$ and $q(X_1, X_2) = X_1 + X_2$, viewing their estimators as elements of $L_1$ or of the space $C_0(\mathbb{R})$ of continuous functions on $\mathbb{R}$ vanishing at infinity. Giné and Mason [13] obtain functional results in $L_p$, locally uniformly in the bandwidth, for general $q(X_1, \ldots, X_m)$. Special cases of the semiparametric time series model considered here have also been studied. Saavedra and Cao [31] consider pointwise convergence of plug-in estimators for the stationary density of moving average processes of order one. Schick and Wefelmeyer [33] obtain asymptotic normality and efficiency and Schick and Wefelmeyer [35] generalize this result to higher-order moving average processes and to functional convergence in $L_1$ and $C_0(\mathbb{R})$; see below for details. Here, we consider general invertible linear processes and obtain $n^{1/2}$-consistency of our estimator for the stationary density in $C_0(\mathbb{R})$.

Specifically, we consider a stationary linear process with infinite-order moving average representation

$$(1.1) \qquad X_t = \varepsilon_t + \sum_{s=1}^{\infty} \varphi_s \varepsilon_{t-s}, \qquad t \in \mathbb{Z},$$

with summable coefficients $\varphi_s$ and independent and identically distributed (i.i.d.) innovations $\varepsilon_t$, $t \in \mathbb{Z}$, having mean zero and finite variance. If the innovations have a density $f$, then $X_0$ has a density, say $h$. The usual estimator of this density from observations $X_1, \ldots, X_n$ of the linear process is a kernel density estimator

$$\tilde{h}(x) = \frac{1}{n} \sum_{j=1}^{n} k_{b_n}(x - X_j), \qquad x \in \mathbb{R},$$

where $k_{b_n} = k(x/b_n)/b_n$ for some kernel $k$ (an integrable function that integrates to 1) and some bandwidth $b_n$ (tending to 0).

Our goal is to construct an $n^{1/2}$-consistent estimator of $h$. For this, we set

$$Y_t = X_t - \varepsilon_t = \sum_{s=1}^{\infty} \varphi_s \varepsilon_{t-s}, \qquad t \in \mathbb{Z}.$$

We must exclude the degenerate case that the observations are i.i.d.:

(C) At least one of the moving average coefficients $\varphi_s$ is nonzero.

$Y_0$ then has a density, say $g$. We have $X_0 = \varepsilon_0 + Y_0$. Since $Y_0$ is independent of $\varepsilon_0$, we can express the density $h$ of $X_0$ as the convolution $h = f * g$ of $f$ and $g$. We obtain an estimator of $h$ as $\hat{h} = \hat{f} * \hat{g}$, where $\hat{f}$ and $\hat{g}$ are estimators of $f$ and $g$. We base these estimators on estimators of the innovations. For this, we require invertibility of the process.

(I) The function $\phi(z) = 1 + \sum_{s=1}^{\infty} \varphi_s z^s$ is bounded and bounded away from zero on the complex unit disk $\{z \in \mathbb{C} : |z| \leq 1\}$.

$\rho(z) = 1/\phi(z) = 1 - \sum_{s=1}^{\infty} \varrho_s z^s$ is then also bounded and bounded away from zero on the complex unit disk. Hence, the innovations have the infinite-order autoregressive representation

$$(1.2) \qquad \varepsilon_t = X_t - \sum_{s=1}^{\infty} \varrho_s X_{t-s}, \qquad t \in \mathbb{Z}.$$

Let $p_n$ be positive integers with $p_n/n \to 0$. For $j = p_n + 1, \ldots, n$, we mimic the innovation $\varepsilon_j$ by the residual

$$\hat{\varepsilon}_j = X_j - \sum_{i=1}^{p_n} \hat{\varrho}_i X_{j-i},$$

where $\hat{\varrho}_i$ is an estimator of $\varrho_i$ for $i = 1, \ldots, p_n$. We then estimate the innovation density by a kernel estimator based on the residuals,

$$\hat{f}(x) = \frac{1}{n - p_n} \sum_{j=p_n+1}^{n} k_{b_n}(x - \hat{\varepsilon}_j), \qquad x \in \mathbb{R},$$

and we estimate the density $g$ by a kernel estimator based on the differences $\hat{Y}_j = X_j - \hat{\varepsilon}_j$,

$$\hat{g}(x) = \frac{1}{n - p_n} \sum_{j=p_n+1}^{n} k_{b_n}(x - \hat{Y}_j), \qquad x \in \mathbb{R}.$$

In addition to (C) and (I), we use the following assumptions:

(Q) the autoregression coefficients satisfy $\sum_{s > p_n} |\varrho_s| = O(n^{-1/2-\zeta})$ for some $\zeta > 0$;

(R) the estimators $\hat{\varrho}_i$ of the autoregression coefficients $\varrho_i$ satisfy

$$\sum_{i=1}^{p_n} (\hat{\varrho}_i - \varrho_i)^2 = O_p(q_n n^{-1})$$

for some $q_n$ with $1 \leq q_n \leq p_n$;

(S) the moving average coefficients satisfy $\sum_{s=1}^{\infty} s|\varphi_s| < \infty$;

(F) the density $f$ has mean zero, a finite fourth moment, is absolutely continuous with a bounded and integrable (almost everywhere) derivative $f'$, and the function $x \mapsto xf'(x)$ is bounded and integrable.

The usual estimators of the autoregression coefficients are the least squares estimators $\hat{\varrho}_1, \ldots, \hat{\varrho}_{p_n}$ which minimize $\sum_{j=p_n+1}^{n}(X_j - \sum_{i=1}^{p_n} \varrho_i X_{j-i})^2$. By Lemma 1, they meet condition (R) with $q_n = p_n$ if, in addition,

$$np_n \sum_{s>p_n} \varrho_s^2 \to 0 \tag{1.3}$$

holds. For smooth parametric models for the autoregression coefficients, we even have (R) with $q_n = 1$, as shown in Section 2.

We denote the number of nonzero coefficients among $\{\varphi_s : s \geq 1\}$ by

$$N = \sum_{s \geq 1} 1[\varphi_s \neq 0].$$

We can then express (C) as $N \geq 1$. If $N$ is finite, then (S) holds and the autocorrelation coefficients decay exponentially. Moreover, (Q) holds with $\zeta = 1$ if $p_n = \log(n)\log(\log n)$.

If we assume that $|\varrho_s| \leq Bs^{-1-\alpha}$ for some $\alpha > 0$, then we have

$$\sum_{s>p_n} |\varrho_s| = O(p_n^{-\alpha}) \quad \text{and} \quad np_n \sum_{s>p_n} \varrho_s^2 = O(np_n^{-2\alpha}).$$

The choice $p_n = n^\beta$ with $2\beta\alpha > 1$ then gives (1.3) and (Q) with $\zeta = \beta\alpha - 1/2$.

Under (C) and (F), the density $h$ is only guaranteed to be twice continuously differentiable. Thus, the optimal rate of nonparametric estimators like the kernel estimator $\tilde{h}$ is $n^{-2/5}$. Our estimator for $h$ is $\hat{h} = \hat{f} * \hat{g}$. We will show that its rate is $n^{-1/2}$. Simulations in [33] for a related estimator in a first-order moving average process show that $\hat{h}$ is better than $\tilde{h}$, even for small sample sizes, and uniformly over a range of bandwidths. We note that our estimator $\hat{h}$ is easy to calculate. Indeed, $\hat{h}(x)$ can be written as the V-statistic

$$\hat{h}(x) = \frac{1}{(n-p_n)^2} \sum_{i=p_n+1}^{n} \sum_{j=p_n+1}^{n} K_{b_n}(x - \hat{\varepsilon}_i - \hat{Y}_j),$$

where $K_b(x) = K(x/b)/b$ and $K = k * k$. Here, we used the fact that $k_b * k_b = K_b$. Thus, it is advantageous to choose a kernel $k$ for which $k * k$ is known.

Smoothness of $g$ and $h$ can be linked to the number $N$. Our main result will thus be formulated in terms of $N$. The following conditions on the kernel and the bandwidth are kept general in order to allow for various smoothness assumptions in terms of an integer $m \geq 2$, where $m - 1$ will play the role of a (known) minimal size for $N$. Under (C), we know that $N \geq 1$, so we can always take $m = 2$.

(B) The sequences $b_n$, $p_n$ and $q_n$ and the exponent $\zeta$ satisfy $p_n q_n b_n^{-1} \times n^{-1/2} \to 0$, $nb_n^{2m} = O(1)$, $n^{1/4} s_n \to 0$ and $n^{1/2} b_n s_n = O(1)$, where $s_n = b_n^{-1/2} n^{-1/2} + p_n q_n b_n^{-5/2} n^{-1} + b_n^{-3/2} n^{-\zeta - 1/2}$.

(K) The kernel $k$ has bounded, continuous and integrable derivatives up to order two and is of type $(m, 2)$, as defined below.

A kernel $k$ is said to be *of type* $(m, c)$ if $\int t^i k(t)\, dt = 0$ for $i = 1, \ldots, m$ and if $\int |t|^{mc} |k(t)|\, dt$ is finite. A kernel satisfying (K) can be chosen to be of the form $p\phi$, where $\phi$ is the standard normal density and $p$ is an appropriate polynomial of degree $m$.

A possible choice of bandwidth is $b_n \sim n^{-1/(2m)}$. Condition (B) is then met if $4m\zeta > 1$ and $p_n q_n n^{-(2m-3)/(4m)} \to 0$ hold. In particular, $p_n = q_n \sim n^\beta$ requires that $8m\beta < 2m - 3$.

Let $\mathbb{G}_n$, $\mathbb{F}_n$ and $\mathbb{H}_n$ denote the processes defined by

$$\mathbb{G}_n(x) = \frac{1}{n - p_n} \sum_{j=p_n+1}^{n} (g(x - \varepsilon_j) - E[g(x - \varepsilon_j)]),$$

$$\mathbb{F}_n(x) = \frac{1}{n - p_n} \sum_{j=p_n+1}^{n} (f(x - Y_j) - E[f(x - Y_j)]),$$

$$\mathbb{H}_n(x) = \sum_{i=1}^{p_n} (\hat{\varrho}_i - \varrho_i) E[X_0 k_{b_n}(x - Y_i)],$$

for $x \in \mathbb{R}$. Let $\|\cdot\|$ denote the supremum norm. We can now state our main result.

THEOREM 1. *Suppose* (I), (Q), (R), (S), (F), (K) *and* (B) *hold. Let* $N \geq m - 1 \geq 1$. *Then*
$$\|\hat{h} - h - \mathbb{F}_n - \mathbb{G}_n + f' * \mathbb{H}_n\| = o_p(n^{-1/2}).$$

The proof is an immediate consequence of the results in Sections 3–10. Write

(1.4) $\quad \hat{h} - h = g * (\hat{f} - f) + f * (\hat{g} - g) + (\hat{f} - f) * (\hat{g} - g).$

Since $f$ is $L_2$-smooth and $g$ is $L_2$-smooth of order $m - 1$, as shown in Section 3, Lemmas 9 and 10 in Section 9 imply $\|\hat{f} - f\|_2 = O_p(s_n) + o(b_n)$, while Lemmas 11 and 12 in Section 10 imply $\|\hat{g} - g\|_2 = O_p(s_n) + o(b_n^{m-1})$. Inequality (4.3) below and condition (B) then give

(1.5) $\quad \|(\hat{f} - f) * (\hat{g} - g)\| \leq \|\hat{f} - f\|_2 \|\hat{g} - g\|_2 = o_p(n^{-1/2}).$

We note that strong consistency of $\hat{f}$ was proved by Robinson [25, 26]. For (finite-order) nonlinear autoregressive models, convergence rates of residual-based kernel estimators were obtained by Liebscher [17] and Müller, Schick and Wefelmeyer [23]. By the smoothness properties of $f$, $g$ and $h$ from Section 3, Theorem 4 in Section 9, applied with $a = g$, gives

(1.6) $\qquad \|g * (\hat{f} - f) - \mathbb{G}_n\| = o_p(n^{-1/2})$

and Theorem 5 in Section 10, applied with $a = f$, gives

$$(1.7) \qquad \|f * (\hat{g} - g) - \mathbb{F}_n + f' * \mathbb{H}_n\| = o_p(n^{-1/2}).$$

Theorem 1 now follows from (1.4)–(1.7).

The sequences $n^{1/2}\mathbb{G}_n$ and $n^{1/2}\mathbb{F}_n$ are tight in $C_0(\mathbb{R})$ by Section 4. Moreover, the sequence $n^{1/2} f' * \mathbb{H}_n$ is tight for the least squares estimators if (1.3) also holds. Indeed, according to Lemma 1 in Section 2, the above assumptions imply that the least squares estimators satisfy

$$(1.8) \qquad \hat{\Delta} = M_n^{-1} \frac{1}{n - p_n} \sum_{j=p_n+1}^{n} \mathbf{X}_{j-1} \varepsilon_j + o_p(n^{-1/2}),$$

where $\hat{\Delta} = (\hat{\varrho}_1 - \varrho_1, \ldots, \hat{\varrho}_{p_n} - \varrho_{p_n})^\top$, $\mathbf{X}_{j-1} = (X_{j-1}, \ldots, X_{j-p_n})^\top$ and $M_n = E[\mathbf{X}_0 \mathbf{X}_0^\top]$. Thus, if (F) holds, then $n^{1/2} f' * \mathbb{H}_n$ is tight in $C_0(\mathbb{R})$ by Theorem 2 in Section 7, applied with $a = f'$. Hence, $n^{1/2}(\hat{h} - h)$ is tight in $C_0(\mathbb{R})$ by the above Theorem 1 and $\hat{h}$ is $n^{1/2}$-consistent in $C_0(\mathbb{R})$. Since the finite-dimensional marginal distributions of $n^{1/2}(\hat{h} - h)$ are asymptotically normal with mean zero, the process $n^{1/2}(\hat{h} - h)$ converges weakly in $C_0(\mathbb{R})$ to a centered Gaussian process with covariance

$$\Gamma(s,t) = \lim_{n \to \infty} \text{Cov}(\mathbb{Z}_n(s), \mathbb{Z}_n(t)), \qquad s, t \in \mathbb{R},$$

where

$$\mathbb{Z}_n(x) = \frac{1}{\sqrt{n}} \sum_{j=1}^{n} (g(x - \varepsilon_j) + f(x - Y_j) - 2h(x) + \varepsilon_j \mathbf{X}_{j-1}^\top M_n^{-1} E[\mathbf{X}_0 f'(x - Y_1)]).$$

We pay a price for $n^{1/2}$-consistency in several respects. One is that we need stronger assumptions on the process, namely invertibility and a sufficiently fast decay of the autoregression coefficients, that is, condition (Q). Another is that we must choose, besides the bandwidth $b_n$, the cut-off index $p_n$. However, our estimator has the advantage that its asymptotic behavior does not depend on $b_n$ and $p_n$, at least in the ranges we allow, while the rate of the usual kernel estimator depends on the bandwidth.

If we strengthen (F) by imposing additional (smoothness) assumptions on $f'$ and use kernels of type $(r, 2)$ for appropriately chosen $r$, the bias terms in the estimation of $f$, $g$ and $h$ can be made smaller, allowing for larger bandwidths and hence weaker assumptions. For example, if $f'$ has bounded variation and a kernel of type $(2m - 1, 2)$ is used, then we can show that $\|f * k_{b_n} - f\|_2 = O(b_n^{3/2})$, $\|g * k_{b_n} - g\|_2 = O(b_n^{2m-5/2})$ and $\|h * k_{b_n} - h\| = O(b_n^{2m-1})$. This allows us to replace the requirements $nb_n^{2m} = O(1)$ and $n^{1/2} b_n s_n = O(1)$ in (B) by $nb_n^{4m-2} \to 0$ and $nb_n^4 = O(1)$. For the choice $b_n = (n \log n)^{1/(4m-2)}$, the requirements of this modified condition (B) are then implied by $p_n q_n (\log n)^{1/2} n^{-(m-1)/(2m-1)} = O(1)$. This allows for larger values of $p_n$ and avoids additional assumptions on $\zeta$.

The paper is organized as follows. In Section 2 we comment more on the assumptions. We also look at the case where we have a parametric model for the autoregressive coefficients and give more details for classical models such as the AR($p$), MA(1) and ARMA(1,1) models. In Section 3 we review expansions in $C_0(\mathbb{R})$ and $L_p$. In Section 4 we give a tightness criterion for sequences of $C_0(\mathbb{R})$-valued random elements and sufficient conditions for tightness of empirical processes based on observations from linear processes. These are used in later sections to show tightness of $n^{1/2}\mathbb{F}_n$, $n^{1/2}\mathbb{G}_n$ and $n^{1/2}f' * \mathbb{H}_n$. An important inequality is established in Section 5. The asymptotic behavior of averages of the form $(n-p_n)^{-1} \sum_{j=p_n+1}^{n} X_{j-i} a_n(x - Y_j)$ and their means is studied in Section 6. Such averages arise in the stochastic expansion of $\hat{g}$. Tightness of $n^{1/2}f' * \mathbb{H}_n$ is established in Section 7. Section 8 shows how well the residuals approximate the true innovations and gives uniform stochastic expansions for residual-based averages of the form $(n-p_n)^{-1} \sum_{j=p_n+1}^{n} a_n(x - \hat{\varepsilon}_j)$ and $(n-p_n)^{-1} \sum_{j=p_n+1}^{n} a_n(x - \hat{Y}_j)$. The kernel estimators $\hat{f}$ and $\hat{g}$ are of this form. In Section 9 we give convergence rates of $\hat{f}$ in $L_2$ and stochastic expansions of functionals $a * \hat{f}$ in $C_0(\mathbb{R})$. Analogous results are given for $\hat{g}$ and $a * \hat{g}$ in Section 10. We have seen above how these results enter the proof of Theorem 1.

**2. Examples.** The following result on the behavior of the least squares estimators is essentially contained in [4].

LEMMA 1. *Assume that* (I), (1.3) *and* $p_n^3/n \to 0$ *hold and that* $f$ *has a finite fourth moment. Then expansion* (1.8) *is valid.*

PROOF. The least squares estimators $(\hat{\varrho}_1, \ldots, \hat{\varrho}_{p_n})^\top$ can be expressed as

$$\hat{M}_n^{-1} \frac{1}{n} \sum_{j=1}^{n} \mathbf{X}_{j-1} X_j \quad \text{with } \hat{M}_n = \frac{1}{n} \sum_{j=1}^{n} \mathbf{X}_{j-1} \mathbf{X}_{j-1}^\top.$$

We can write the error term in (1.8) as $(\hat{M}_n^{-1} - M_n^{-1}) A_n - \hat{M}_n^{-1} B_n$ with

$$A_n = \frac{1}{n} \sum_{j=1}^{n} \mathbf{X}_{j-1} \varepsilon_j \quad \text{and} \quad B_n = \frac{1}{n} \sum_{j=1}^{n} \mathbf{X}_{j-1} \sum_{i>p_n} \varrho_i X_{j-i}.$$

By (2.13) of Berk [4],

$$E[|B_n|^2] = O\left(p_n \sum_{i>p_n} \varrho_i^2\right)$$

and by the relation immediately preceding his (2.17), we have $E[|A_n|^2] = O(p_n n^{-1})$. By his Lemma 3, we have $p_n^{1/2} \|\hat{M}_n^{-1} - M_n^{-1}\|_* = o_p(1)$, where $\|M\|_* = \sup_{|x| \leq 1} |Mx|$ is the operator norm of a matrix $M$. By his (2.14), both $\|M_n\|_*$ and $\|M_n^{-1}\|_*$ are bounded. Combining the above, we obtain

$$(\hat{M}_n^{-1} - M_n^{-1}) A_n = o_p(p_n^{-1/2}) O_p(p_n^{1/2} n^{-1/2}) = o_p(n^{-1/2}),$$

$$\hat{M}_n^{-1} B_n = O_p\left(p_n^{1/2} \left(\sum_{i > p_n} \varrho_i^2\right)^{1/2}\right) = o_p(n^{-1/2}).$$

The result follows. □

Of special interest is the case where we have a parametric model for the autocorrelation coefficients, that is, there are functions $r_1, r_2, \ldots$ from an open subset $\Theta$ of $\mathbb{R}^q$ into $\mathbb{R}$ such that $\varrho_i = r_i(\vartheta)$ for all $i$ and some unknown $\vartheta$ in $\Theta$. We can then take $\hat{\varrho}_i = r_i(\hat{\vartheta})$ for all $i$ and some estimator $\hat{\vartheta}$ of $\vartheta$. Now, let us impose the following conditions:

(R1) the estimator $\hat{\vartheta}$ of $\vartheta$ is $n^{1/2}$-consistent, that is, $\hat{\vartheta} - \vartheta = O_p(n^{-1/2})$;

(R2) the functions $r_1, r_2, \ldots$ are differentiable at $\vartheta$ with gradients $\dot{r}_1(\vartheta)$, $\dot{r}_2(\vartheta), \ldots$ and

$$\sum_{i=1}^{\infty} (r_i(\vartheta + s) - r_i(\vartheta) - \dot{r}_i(\vartheta)^\top s)^2 = o(|s|^2) \quad \text{and} \quad \sum_{i=1}^{\infty} |\dot{r}_i(\vartheta)|^2 < \infty.$$

These conditions imply (R) with $q_n = 1$. If (C) and (F) are also met, one obtains (see Theorem 3 in Section 7) that

$$\|f' * \mathbb{H}_n - (\hat{\vartheta} - \vartheta)^\top \Lambda\| = o_p(n^{-1/2})$$

with

$$\Lambda(x) = \sum_{i=1}^{\infty} \dot{r}_i(\vartheta) E[X_0 f'(x - Y_i)], \qquad x \in \mathbb{R}.$$

Thus, if (I), (Q), (R1), (R2), (S), (F), (K), (B) and $N \geq m - 1$ hold, we have the expansion

(2.1) $$\|\hat{h} - h - \mathbb{F}_n - \mathbb{G}_n + (\hat{\vartheta} - \vartheta)^\top \Lambda\| = o_p(n^{-1/2})$$

and tightness of $n^{1/2}(\hat{h} - h)$. Weak convergence of $n^{1/2}(\hat{h} - h)$ in $C_0(\mathbb{R})$ can now be established under mild additional assumptions on $\hat{\vartheta}$.

Let us now look at three special cases, namely AR($p$), MA(1) and ARMA(1,1). In these examples, the moving average and autoregression coefficients decay exponentially, so (S) holds and the choice $p_n \sim \log(n) \log(\log(n))$ guarantees (Q) with $\zeta = 1$. We can then take $m = 2$ and $b_n \sim n^{-1/4}$.

EXAMPLE 1. Let $X_t = \vartheta_1 X_{t-1} + \cdots + \vartheta_p X_{t-p} + \varepsilon_t$ be an AR($p$) process with $\vartheta_p \neq 0$ and such that the polynomial $\varrho(z) = 1 - \sum_{i=1}^{p} \vartheta_i z^i$ has no roots in the (complex) unit disk. Set $\vartheta = (\vartheta_1, \ldots, \vartheta_p)^\top$ and $\tilde{X}_{t-1} = (X_{t-j}, \ldots, X_{t-p})^\top$. We can then write the model as $X_t = \vartheta^\top \tilde{X}_{t-1} + \varepsilon_t$. The representation (1.2) holds with $\varrho_s = r_s(\vartheta) = \vartheta_s$ for $s \leq p$ and $\varrho_s = r_s(\vartheta) = 0$ for $s > p$. By our assumptions on $\varrho(z)$, the moving average representation (1.1) holds with $\varphi_s$ being the coefficients of $1/\varrho(z) = \sum_{s=1}^{\infty} \varphi_s z^k$ and $Y_t = X_t - \varepsilon_t = \vartheta^\top \tilde{X}_{t-1}$. Since $\vartheta = 0$ is ruled out, we have (C). Moreover, the moving average coefficients decay exponentially, implying (S). Let $\hat\vartheta$ be an $n^{1/2}$-consistent estimator of $\vartheta$. We estimate the innovations $\varepsilon_j$ by the residuals $\hat\varepsilon_j = X_j - \hat\vartheta^\top \tilde{X}_{j-1}$. Here, (R2) holds with $\dot{r}_i(\vartheta) = e_i$, the $i$th unit vector, for $i \leq p$ and with $\dot{r}_i(\vartheta) = 0$ for $i > p$, and we find $\Lambda(x) = E[\tilde{X}_0 f'(x - \vartheta^\top \tilde{X}_0)]$. A simple estimator for $\vartheta$ is the least squares estimator

$$\hat\vartheta = \left( \sum_{j=p+1}^{n} \tilde{X}_{j-1} \tilde{X}_{j-1}^\top \right)^{-1} \sum_{j=p+1}^{n} \tilde{X}_{j-1} X_j.$$

With $M = E[\tilde{X}_0 \tilde{X}_0^\top]$, $\hat\vartheta$ has the stochastic expansion

$$\hat\vartheta = \vartheta + M^{-1} \frac{1}{n} \sum_{j=1}^{n} \tilde{X}_{j-1} \varepsilon_j + o_p(n^{-1/2}).$$

With this choice of $\hat\vartheta$, we obtain, in particular, that $n^{1/2}(\hat h - h)$ converges weakly in $C_0(\mathbb{R})$ to a centered Gaussian process. In this example, we can take $p_n = p$.

EXAMPLE 2. Let $X_t = \varepsilon_t + \vartheta \varepsilon_{t-1}$ be an MA(1) process with $|\vartheta| < 1$ and $\vartheta \neq 0$. The moving average representation (1.1) then holds with $\varphi_1 = \vartheta$ and $\varphi_s = 0$ for $s > 1$, and (C) holds, as $\vartheta \neq 0$. The representation (1.2) holds with $\varrho_s = r_s(\vartheta) = -(-\vartheta)^s$. Let $\hat\vartheta$ be an $n^{1/2}$-consistent estimator of $\vartheta$. We estimate the innovations $\varepsilon_j$ by the residuals $\hat\varepsilon_j = X_j + \sum_{i=1}^{p_n} (-\hat\vartheta)^i X_{j-i}$. It is easy to check that (R2) holds with $\dot r_s(\vartheta) = s(-\vartheta)^{s-1}$. We have $Y_t = X_t - \varepsilon_t = \vartheta \varepsilon_{t-1}$ and, therefore, $E[X_0 f'(x - Y_i)] = 0$ for $i > 1$. Thus, the expansion (2.1) holds with $\Lambda(x) = E[X_0 f'(x - Y_1)] = E[\varepsilon_0 f'(x - \vartheta\varepsilon_0)]$. In particular, if $\hat\vartheta$ is asymptotically linear, then $n^{1/2}(\hat h - h)$ converges weakly in $C_0(\mathbb{R})$ to a centered Gaussian process. Our estimator $\hat h$ is asymptotically equivalent to the estimator

$$\hat h_{SC}(x) = \int \hat f(x - \hat\vartheta y) \hat f(y) \, dy$$

considered by Saavedra and Cao [31]. This estimator can be written

$$\hat h_{SC}(x) = \frac{1}{n^2 b_n} \sum_{i=1}^{n} \sum_{j=1}^{n} L_{\hat\vartheta}\left( \frac{x - \varepsilon_i - \hat\vartheta \varepsilon_j}{b_n} \right)$$

with $L_\vartheta(x) = \int k(x - \vartheta y)k(y)\,dy$. The kernel $L_{\hat\vartheta}$ can be replaced by a general (nonrandom) kernel $k$. The U-statistic version of the resulting estimator,

$$\hat{h}_{SW} = \sum_{\substack{i,j=1 \\ i \neq j}}^{n} k_{b_n}(x - \varepsilon_i - \hat\vartheta\varepsilon_j),$$

is studied in [33], where a pointwise version of the above stochastic expansion is proved. Schick and Wefelmeyer [35] generalize the result to MA($q$) and show that the expansion holds uniformly and in $L_1$.

EXAMPLE 3. Let $X_t = \alpha X_{t-1} + \varepsilon_t + \beta\varepsilon_{t-1}$ be an ARMA(1,1) process with $|\alpha|, |\beta| < 1$ and $\alpha + \beta \neq 0$. The moving average representation (1.1) then holds with $\varphi_s = (\alpha + \beta)\alpha^{s-1}$ and the autoregressive representation (1.2) holds with $\varrho_s = r_s(\alpha,\beta) = (\alpha + \beta)(-\beta)^{s-1}$. The requirement that $\alpha + \beta \neq 0$ gives $\varphi_1 \neq 0$ and, therefore, (C). We have $Y_t = X_t - \varepsilon_t = \sum_{s=1}^{\infty}(\alpha + \beta)\alpha^{s-1}\varepsilon_{t-s}$. Let $\hat\alpha$ and $\hat\beta$ be $n^{1/2}$-consistent estimators of $\alpha$ and $\beta$, respectively. We estimate the innovations $\varepsilon_j$ by the residuals

$$\hat\varepsilon_j = X_j - (\hat\alpha + \hat\beta)\sum_{i=1}^{p_n}(-\hat\beta)^{i-1}X_{j-i}.$$

Here, (R2) holds with $\dot{r}_s(\alpha,\beta) = ((-\beta)^{s-1}, -(s-1)\alpha(-\beta)^{s-2} + s(-\beta)^{s-1})^\top$. Thus, the expansion (2.1) holds with $\hat\vartheta = (\hat\alpha,\hat\beta)^\top$ and

$$\Lambda(x) = \sum_{s=1}^{\infty}\begin{pmatrix} (-\beta)^{s-1} \\ -(s-1)\alpha(-\beta)^{s-2} + s(-\beta)^{s-1} \end{pmatrix} E[X_0 f'(x - Y_s)].$$

In particular, if $\hat\alpha$ and $\hat\beta$ are asymptotically linear, then $n^{1/2}(\hat{h} - h)$ converges weakly in $C_0(\mathbb{R})$ to a centered Gaussian process.

**3. Smoothness.** Here, we shall address smoothness of $f$, $g$ and $h = f * g$. For this, we assume that $N \geq r$ for some positive integer $r$. We can then express $Y_0 = \sum_{i=1}^{r}\varphi_{\tau_i}\varepsilon_{-\tau_i} + Z$, where $\tau_1, \ldots, \tau_r$ are the indices of the first $r$ nonzero terms among $\{\varphi_s : s \geq 1\}$ and $Z = \sum_{s > \tau_r}\varphi_s\varepsilon_{-s}$. For $t \neq 0$, define densities $f_t$ and $\bar{f}_t$ by $f_t(x) = f(x/t)/|t|$ and $\bar{f}_t(x) = E[f_t(x - Z)]$. Since the innovations are independent with density $f$, we find that the density $g$ of $Y_0$ equals $\bar{f}_{\tau_1}$ if $r = 1$ and equals the convolution $f_{\tau_1} * \cdots * f_{\tau_{r-1}} * \bar{f}_{\tau_r}$ if $r > 1$.

Let $\mathcal{A}$ denote the class of absolutely continuous functions with a bounded and integrable almost everywhere derivative. Let $\mathcal{A}_p$ denote the class of absolutely continuous functions with an almost everywhere derivative in $L_p$, $p \in [1, \infty)$. It follows from (F) that $f$ belongs to $\mathcal{A}$ and, hence, to $\mathcal{A}_p$ for each $p \in [1, \infty)$. Elements of $\mathcal{A}$ are Lipschitz, while elements $a$ of $\mathcal{A}_p$ are $L_p$-Lipschitz with constant $C = \|a'\|_p$, that is,

$$\|a(\cdot - t) - a\|_p \leq C|t|, \qquad t \in \mathbb{R}.$$

Indeed, we can express

$$a(x+t) - a(x) = t \int_0^1 a'(x+st)\,ds$$

and thus obtain from Jensen's inequality and Fubini's theorem that

$$\int |a(x+t) - a(x)|^p\,dx \leq |t|^p \int_0^1 \int |a'(x+st)|^p\,dx\,ds = |t|^p \|a'\|_p^p, \qquad t \in \mathbb{R}.$$

A more careful analysis shows that elements $a$ of $\mathcal{A}_p$ are $L_p$-smooth,

$$\|a(\cdot - t) - a + ta'\|_p \leq |t| w_{p,a'}(|t|), \qquad t \in \mathbb{R}.$$

Here, $w_{p,v}$ denotes the $L_p$-modulus of continuity of a measurable function $v$, defined by

$$w_{p,v}(\delta) = \sup_{|t| \leq \delta} \|v(\cdot - t) - v\|_p, \qquad \delta \geq 0.$$

If $v$ belongs to $L_p$, then $w_{p,v}$ is bounded by $2\|v\|_p$ and $w_{p,v}(\delta) \to 0$ as $\delta \to 0$, by the translation continuity in $L_p$, for which we refer to Theorem 9.5 in [30]. Also, recall that the *modulus of continuity* of a function $v$ is defined by

$$w_v(\delta) = \sup_{x,y \in \mathbb{R}, |y-x| \leq \delta} |v(y) - v(x)| \leq \sup_{|t| \leq \delta} \|v(\cdot - t) - v\|, \qquad \delta \geq 0.$$

If $v$ belongs to $C_0(\mathbb{R})$, then $w_v$ is bounded by $2\|v\|$ and $w_v(\delta) \to 0$ as $\delta \to 0$.

Assume now that $f$ belongs to $\mathcal{A}$. Then the densities $f_t$ and $\bar{f}_t$ for $t \neq 0$ will also belong to $\mathcal{A}$. This immediately gives that $g$ belongs to $\mathcal{A}$ if $r = 1$. Hence, $g$ is $L_p$-smooth for each $1 \leq p < \infty$. Now, assume that $r > 1$. Set $g_i = f'_{\tau_1} * \cdots * f'_{\tau_i} * f_{\tau_{i+1}} * \cdots * f_{\tau_{r-1}} * \bar{f}_{\tau_r}$ for $i = 1, \ldots, r-1$ and $g_r = f'_{\tau_1} * \cdots * f'_{\tau_{r-1}} * \bar{f}'_{\tau_r}$. These functions are integrable, bounded and uniformly continuous. The last two properties stem from the fact that the convolution of a bounded function $u$ with an integrable function $v$ is bounded and uniformly continuous in view of the bounds $\|u*v\| \leq \|u\|\|v\|_1$ and $w_{u*v}(\delta) \leq \|u\|w_{1,v}(\delta)$. It is now easy to check that $g_i$ is the $i$th derivative of $g$. Thus, we have the identity

$$g(x+t) - g(x) - \sum_{i=1}^r \frac{t^i}{i!} g_i(x) = \frac{t^r}{r!} \int_0^1 (g_r(x+st) - g_r(x)) r(1-s)^{r-1}\,ds.$$

Since $g_r$ belongs to $L_p$, we obtain from Jensen's inequality and Fubini's theorem, as above, that

$$(3.1) \qquad \left\| g(\cdot + t) - g - \sum_{i=1}^r \frac{t^i}{i!} g_i \right\|_p \leq \frac{|t|^r}{r!} w_{p,g_r}(|t|), \qquad t \in \mathbb{R}.$$

If (3.1) holds and $g_r \in L_p$, then we say that $g$ is $L_p$-smooth of order $r$. This property reduces to $L_p$-smoothness if $r = 1$.

Since $h$ equals $f * g$, the above arguments show that $h$ is $(r+1)$-times continuously differentiable with bounded, integrable and uniformly continuous derivatives. This implies that

$$(3.2) \qquad \left\| h(\cdot + t) - h - \sum_{i=1}^{r+1} \frac{t^i}{i!} h^{(i)} \right\| \leq \frac{|t|^{r+1}}{(r+1)!} w_{h^{(r+1)}}(|t|), \qquad t \in \mathbb{R}.$$

If (3.2) holds and $h^{(r+1)}$ is bounded and uniformly continuous, we say that $h$ is *smooth of order* $r+1$.

Let us now summarize our findings.

COROLLARY 1. *Let $f$ belong to $\mathcal{A}$ and $N \geq r \geq 1$. Then $f$ is $L_2$-smooth, $g$ belongs to $\mathcal{A}$ and is $L_2$-smooth of order $r$ and $h$ is smooth of order $r+1$.*

COROLLARY 2. *Let $a$ be $L_2$-smooth of order $r$ and let $k$ be a kernel of type $(m, 2)$ with $m \geq r$. Then $\|a * k_{b_n} - a\|_2 = o(b_n^r)$.*

COROLLARY 3. *Let $a$ be smooth of order $r$ and let $k$ be a kernel of type $(m, 1)$ with $m \geq r$. Then $\|a * k_{b_n} - a\| = o(b_n^r)$.*

**4. Weak convergence in $C_0(\mathbb{R})$.** In this section, we address weak convergence of sequences of random elements in the space $C_0(\mathbb{R})$ of continuous functions vanishing at (plus and minus) infinity, endowed with the supremum norm $\|\cdot\|$. To establish tightness, we use the following characterization of compact subsets of $C_0(\mathbb{R})$.

LEMMA 2. *A closed subset $A$ of $C_0(\mathbb{R})$ is compact if and only if*

$$\limsup_{\delta \downarrow 0} \sup_{a \in A} \sup_{|z-y| \leq \delta} |a(z) - a(y)| = 0,$$

$$\lim_{K \to \infty} \sup_{a \in A} \sup_{|z| \geq K} |a(z)| = 0.$$

A proof of this lemma is given in [34]. From the lemma, we immediately obtain the following characterization of tightness.

COROLLARY 4. *A sequence $\mathbb{A}_n$ of $C_0(\mathbb{R})$-valued random elements is tight if and only if for every $\varepsilon > 0$ and $\eta > 0$, there exist a $\delta > 0$ and a $K < \infty$ such that*

$$(4.1) \qquad \sup_n P\left( \sup_{|z-y| \leq \delta} |\mathbb{A}_n(z) - \mathbb{A}_n(y)| > \varepsilon \right) < \eta,$$

$$(4.2) \qquad \sup_n P\left( \sup_{|z| \geq K} |\mathbb{A}_n(z)| > \varepsilon \right) < \eta.$$

Once tightness is established, weak convergence follows from the convergence of the finite-dimensional distributions.

Let $a_1$ and $a_2$ be two square-integrable functions. Then $a_1 * a_2$ belongs to $C_0(\mathbb{R})$. Indeed, an application of the Cauchy–Schwarz inequality and a substitution yield

$$\|a_1 * a_2\| \leq \|a_1\|_2 \|a_2\|_2. \tag{4.3}$$

Hence, $a_1 * a_2$ is bounded. Furthermore,

$$\|a_1 * a_2(\cdot - t) - a_1 * a_2\| \leq \|a_1(\cdot - t) - a_1\|_2 \|a_2\|_2. \tag{4.4}$$

Since $a_1$ is square-integrable, we obtain from the translation continuity of square-integrable functions (see, e.g., [30], Theorem 9.5) that $\|a_1(\cdot - t) - a_1\|_2 \to 0$ as $t \to 0$. This shows that $a_1 * a_2$ is uniformly continuous. Finally, write $\chi_K(y) = \mathbf{1}[|y| > K]$ and $a_1 * a_2 = a_1 * (a_2(1 - \chi_K)) + a_1 * (a_2 \chi_K)$. Since $|x - y| > K$ if $|x| > 2K$ and $|y| \leq K$, we obtain

$$\sup_{|x|>2K} |a_1 * a_2(x)| \leq \|a_1 \chi_K\|_2 \|a_2\|_2 + \|a_1\|_2 \|a_2 \chi_K\|_2. \tag{4.5}$$

Hence $a_1 * a_2$ vanishes at infinity. The above shows that $a_1 * a_2$ is in $C_0(\mathbb{R})$.

If $a$ is a square-integrable function and $\mathbb{D}_n$ is a sequence of $L_2$-valued random elements, then inequalities (4.3)–(4.5) yield

$$\|a * \mathbb{D}_n(\cdot - t) - a * \mathbb{D}_n\| \leq \|a(\cdot - t) - a\|_2 \|\mathbb{D}_n\|_2,$$

$$\sup_{|x|>2K} |a * \mathbb{D}_n(x)| \leq \|a \chi_K\|_2 \|\mathbb{D}_n\|_2 + \|a\|_2 \|\mathbb{D}_n \chi_K\|_2.$$

This shows that the $C_0(\mathbb{R})$-valued sequence $a * \mathbb{D}_n$ is tight if $\|\mathbb{D}_n\|_2 = O_p(1)$ and if for all positive $\varepsilon$ and $\eta$, there exists a $K$ such that $\sup_n P(\|\mathbb{D}_n \chi_K\|_2 > \varepsilon) < \eta$. In view of the Markov inequality, a sufficient condition for these two statements is the following condition.

(T) There exists an integrable $\Psi$ such that $E[\mathbb{D}_n^2(x)] \leq \Psi(x)$ for all $x \in \mathbb{R}$.

Now, let $\xi_1, \xi_2, \ldots$ be a stationary sequence of random variables with distribution function $D$ and let

$$\mathbb{D}_n(x) = n^{-1/2} \sum_{j=1}^n (\mathbf{1}[\xi_j \leq x] - D(x)), \qquad x \in \mathbb{R},$$

be the associated empirical process. If $A$ is absolutely continuous with an almost everywhere derivative $A'$ that is both integrable and square-integrable, then we can express

$$\mathbb{A}_n(x) = n^{-1/2} \sum_{j=1}^n (A(x - \xi_j) - E[A(x - \xi_j)]) = \int A(x - y) \, d\mathbb{D}_n(y)$$

as

$$\mathbb{A}_n(x) = \int A'(x-y)\mathbb{D}_n(y)\,dy = A' * \mathbb{D}_n(x), \qquad x \in \mathbb{R}.$$

Thus, the sequence $\mathbb{A}_n$ will be tight if we can show that condition (T) holds. In the following, we give sufficient conditions for (T).

(a) If $\xi_1, \xi_2, \ldots$ are independent, then condition (T) holds if the random variables have a finite mean. Indeed, we have the identity $E[\mathbb{D}_n^2(x)] = D(x)(1 - D(x))$ and $D(1-D)$ is integrable if and only if the $\xi_j$ have finite mean.

(b) Now assume that $\xi_1, \xi_2, \ldots$ come from a linear process

$$\xi_t = \sum_{s=0}^{\infty} d_s U_{t-s}, \qquad t \in \mathbb{Z},$$

where the innovations $U_t$, $t \in \mathbb{Z}$, are i.i.d. with finite mean, the coefficients $d_0, d_1, \ldots$ are summable and $d_0 \neq 0$. Then condition (T) holds if $\sum_{s=0}^{\infty}(1+s)|d_s| < \infty$. This follows from Corollary 7.1 in [36].

**5. A bound.** Let $U_t$, $t \in \mathbb{Z}$, be independent and identically distributed random variables with finite mean. For summable coefficients $c_0, c_1, \ldots$ and $d_0, d_1, \ldots$ with $d_0 \neq 0$, let us consider the linear processes

$$S_t = \sum_{s=0}^{\infty} c_s U_{t-s} \quad \text{and} \quad T_t = \sum_{s=0}^{\infty} d_s U_{t-s}, \qquad t \in \mathbb{Z}.$$

For a measurable function $a$, we define

$$K(x) = n^{-1/2} \sum_{j=1}^{n}(a(x - T_j) - E[a(x - T_j)]),$$

$$H(x) = n^{-1/2} \sum_{j=1}^{n}(S_j a(x - T_j) - E[S_j a(x - T_j)]), \qquad x \in \mathbb{R}.$$

Let $U = U_0$ and set

$$\alpha = \sum_{j=0}^{\infty} |c_j| \quad \text{and} \quad D = \sum_{j=0}^{\infty}(j+1)|d_j| = \sum_{j=0}^{\infty}\sum_{s=j}^{\infty}|d_s|.$$

In their Lemma 7.3, Schick and Wefelmeyer [36] show the following result.

LEMMA 3. *Suppose $a$ is bounded and $L_1$-Lipschitz with constant $L$. Let $D$ be finite. Then*

$$\int E[K^2(x)]\,dx \leq 4L\|a\|DE[|U|].$$

We shall now obtain a similar result for the process $H$.

LEMMA 4. *Suppose $a$ is bounded and $L_1$-Lipschitz with constant $L$ and $U$ has a finite second moment. Let $D$ be finite. Then*

$$\int E[H^2(x)]\,dx \le 8L\|a\|\alpha^2 D E[|U|]E[U^2].$$

PROOF. We can write $H(x) = n^{-1/2}\sum_{j=1}^n(Z_j(x) - E[Z_j(x)])$, where

$$Z_j(x) = S_j a(x - T_j), \qquad x \in \mathbb{R}.$$

Now, set

$$S_j^* = \sum_{s=0}^{j-1} c_s U_{j-s}, \qquad \bar{S}_j = \sum_{s=j}^{\infty} c_s U_{j-s},$$

$$T_j^* = \sum_{s=0}^{j-1} d_s U_{j-s}, \qquad \bar{T}_j = \sum_{s=j}^{\infty} d_s U_{j-s}.$$

We can then write

$$Z_j(x) = S_j^* a(x - T_j^* - \bar{T}_j) + \bar{S}_j a(x - T_j^* - \bar{T}_j)$$

and obtain, with $\mathcal{F}$ denoting the $\sigma$-field generated by $\{U_t : t \le 0\}$, that

(5.1) $$\bar{Z}_j(x) = E(Z_j(x)|\mathcal{F}) = a_j^*(x - \bar{T}_j) + \bar{S}_j a_j(x - \bar{T}_j),$$

where $a_j^*$ and $a_j$ are the functions defined by

$$a_j^*(x) = E[S_j^* a(x - T_j^*)] \quad \text{and} \quad a_j = E[a(x - T_j^*)], \qquad x \in \mathbb{R}.$$

These functions inherit the $L_1$-Lipschitz property of $a$. More precisely, we have the bounds

(5.2) $$\|a_j^*(\cdot - t) - a_j^*\|_1 \le E[|S_j^*|]L|t| \le BL|t| \quad \text{and}$$
$$\|a_j(\cdot - t) - a_j\|_1 \le L|t|,$$

where $B = \alpha E[|U|]$. To simplify notation, we abbreviate $S_0$ by $S$, $T_0$ by $T$ and $Z_0$ by $Z$. Using stationarity and a conditioning argument, we obtain

$$E[H^2(x)] = \text{Var}(Z(x)) + \frac{2}{n}\sum_{j=1}^{n-1}(n-j)\text{Cov}(Z(x),\bar{Z}_j(x)) \le 2\sum_{j=0}^{\infty}\Gamma_j(x),$$

where, in view of (5.1), $\Gamma_j(x)$ can be taken to be

$$\Gamma_j(x) = E[|Z(x) - E[Z(x)]||a_j^*(x - \bar{T}_j) - a_j^*(x) + \bar{S}_j(a_j(x - \bar{T}_j) - a_j(x))|].$$

Since $a$ is bounded, we derive the bounds $|Z(x)| \leq |S|\|a\|$ and $|E[Z(x)]| \leq E[|S|]\|a\|$ for $x \in \mathbb{R}$. This, $E[|S|] \leq B = \alpha E[|U|]$ and (5.2) yield that

$$\|\Gamma_j\|_1 \leq \|a\| E[(|S| + E[|S|])(BL|\bar{T}_j| + LE[|\bar{S}_j \bar{T}_j|])]$$

$$\leq \|a\| BL \left( \sum_{s \geq 0} |d_{s+j}| E[(|S| + E[|S|])|U_{-s}|] + 2 \sum_{s,t \geq j} |c_t||d_s| E[U^2] \right)$$

$$\leq \|a\| BL(2\alpha E[U^2] + 2\alpha E[U^2]) \sum_{s \geq j} |d_s|.$$

In view of $B = \alpha E[|U|]$ and the definition of $D$, the desired result is now immediate. $\square$

**6. An auxiliary result.** Let $X_t$ be a linear process as in (1.1). Let $a_n$ be an integrable function that belongs to $\mathcal{A}_1$. For $i = 1, 2, \ldots$, set

$$\hat{a}_{n,i}(x) = \frac{1}{n - p_n} \sum_{j=p_n+1}^{n} X_{j-i} a_n(x - Y_j), \qquad x \in \mathbb{R},$$

$$\bar{a}_{n,i}(x) = E[\hat{a}_{n,i}(x)] = E[X_0 a_n(x - Y_i)], \qquad x \in \mathbb{R}.$$

In this section, we study the behavior of $\hat{a}_{n,i}$ and its expectation $\bar{a}_{n,i}$ in $L_2$. The results developed here will be used in later sections with $a_n = k_{b_n}$ or $a_n = k'_{b_n}$.

From Lemma 4, we immediately obtain the following result.

LEMMA 5. *Suppose* (C) *and* (S) *hold. Then there exists a finite constant $A$ such that*

$$\int \operatorname{Var}(\hat{a}_{n,i}(x)) \, dx \leq A \|a_n\| \|a'_n\|_1, \qquad i = 1, 2, \ldots.$$

We denote the index of the first nonzero moving average coefficient by

$$\tau = \inf\{s \geq 1 : \varphi_s \neq 0\}.$$

Under (C), $\tau$ is finite. Let $Z_j = Y_j - \varphi_\tau \varepsilon_{j-\tau}$. A conditioning argument shows that

$$\bar{a}_{n,i}(x) = \mathbf{1}[i = \tau] E[v_n(x - Z_i)] + E[X_0 u_n(x - Z_i)]$$

with

$$u_n(x) = E[a_n(x - \varphi_\tau \varepsilon_0)] \quad \text{and} \quad v_n(x) = E[\varepsilon_0 a_n(x - \varphi_\tau \varepsilon_0)], \qquad x \in \mathbb{R}.$$

Then $u_n = a_n * \psi_0$ and $v_n = a_n * \psi_1$, where

$$(6.1) \quad \psi_0(x) = \frac{1}{|\varphi_\tau|} f\left(\frac{x}{\varphi_\tau}\right) \quad \text{and} \quad \psi_1(x) = \frac{1}{|\varphi_\tau|} \frac{x}{\varphi_\tau} f\left(\frac{x}{\varphi_\tau}\right), \qquad x \in \mathbb{R}.$$

Under assumption (F), $\psi_0$ and $\psi_1$ belong to $\mathcal{A}$.

If $u_n$ converges in $L_2$ to some $u$ and $v_n$ converges in $L_2$ to some $v$, then we find that $\bar{a}_{n,i}$ converges in $L_2$ to $\bar{a}_i$, where

$$\bar{a}_i(x) = \mathbf{1}[i=\tau]E[v(x-Z_i)] + E[X_0 u(x-Z_i)], \qquad x \in \mathbb{R}.$$

Actually, a stronger statement is possible.

LEMMA 6. *Let* (C), (S) *and* (F) *hold. Suppose that there exist square-integrable functions* $u$ *and* $v$ *with* $u$ *in* $\mathcal{A}_2$ *such that* $\|a_n * \psi_1 - v\|_2 \to 0$, $\|a_n * \psi_0 - u\|_2 \to 0$ *and* $\|a_n * \psi_0' - u'\|_2 \to 0$. *Then*

$$\sum_{i=1}^{\infty} \|\bar{a}_{n,i} - \bar{a}_i\|_2^2 \to 0 \quad and \quad \sum_{i=1}^{\infty} \|\bar{a}_i\|_2^2 < \infty.$$

PROOF. For $i > \tau$ and $w \in \mathcal{A}_2$, we have

$$E[X_0 w(x - Z_i)] = E[X_0 (w(x-Z_i) - w(x-\bar{Z}_i))]$$

with $\bar{Z}_i = \sum_{\tau < s < i} \varphi_s \varepsilon_{i-s}$ and, hence,

$$\int (E[X_0 w(x-Z_i)])^2 \, dx \le E[X_0^2] \int E[(w(x-Z_i) - w(x-\bar{Z}_i))^2] \, dx$$

$$\le E[X_0^2] \|w'\|_2^2 E[(Z_i - \bar{Z}_i)^2]$$

$$= E[X_0^2] \|w'\|_2^2 E[\varepsilon_0^2] \sum_{s=i}^{\infty} \varphi_s^2.$$

With $w = a_n * \psi_0 - u$ and assumption (S), we obtain

$$\sum_{i > \tau} \|\bar{a}_{n,i} - \bar{a}_i\|_2^2 \le E[X_0^2] E[\varepsilon_0^2] \|a_n * \psi_0' - u'\|_2^2 \sum_{s > \tau} s \varphi_s^2 \to 0,$$

and with $w = u$, we obtain

$$\sum_{i > \tau} \|\bar{a}_i\|_2^2 \le E[X_0^2] E[\varepsilon_0^2] \|u'\|_2^2 \sum_{s > \tau} s \varphi_s^2 < \infty.$$

The desired results are now immediate, as $\bar{a}_{n,i}$ converges in $L_2$ to $\bar{a}_i$ for $i \le \tau$. □

REMARK 1. The assumptions on $a_n$ of the previous lemma hold with $u = a * \psi_0$ and $v = a * \psi_1$ if $a_n$ converges in $L_2$ to some $a$. They hold with $u = \psi_0$ and $v = \psi_1$ if $a_n = k_{b_n}$. In the first case, $\bar{a}_i = a * \delta_i$, and in the second case, $\bar{a}_i = \delta_i$, where

(6.2) $\qquad \delta_i(x) = \mathbf{1}[i=\tau]E[\psi_1(x-Z_0)] + E[X_0 \psi_0(x-Z_i)].$

**7. Tightness of $n^{1/2}a * \mathbb{H}_n$.** Let us now address tightness of $n^{1/2}a * \mathbb{H}_n$ for some square-integrable $a$. For such an $a$, we have, with $a_n = a * k_{b_n}$,

$$a * \mathbb{H}_n(x) = \sum_{i=1}^{p_n} (\hat{\varrho}_i - \varrho_i) E[X_0 a_n(x - Y_i)] = \hat{\Delta}^\top E[\mathbf{X}_0 a_n(x - Y_1)], \qquad x \in \mathbb{R}.$$

Recall that $\hat{\Delta} = (\hat{\varrho}_1 - \varrho_1, \ldots, \hat{\varrho}_{p_n} - \varrho_{p_n})^\top$ and $\mathbf{X}_{j-1} = (X_{j-1}, \ldots, X_{j-p_n})^\top$. We shall first treat the case where (1.8) holds. As seen in the proof of Lemma 1, the dispersion matrix $M_n = E[\mathbf{X}_0 \mathbf{X}_0^\top]$ is invertible and the operator norm of its inverse $M_n^{-1}$ is bounded. Hence, there exists a constant $K$ such that for all $n$,

$$(7.1) \quad c_n^\top M_n c_n \leq K |c_n|^2 \quad \text{and} \quad c_n^\top M_n^{-1} c_n \leq K |c_n|^2, \qquad c_n \in \mathbb{R}^{p_n}.$$

Let $\boldsymbol{\delta} = (\delta_1, \ldots, \delta_{p_n})^\top$ with $\delta_i$ as defined in (6.2). Now, set

$$\mathbb{J}_n(x) = \frac{1}{n - p_n} \sum_{j=p_n+1}^{n} \varepsilon_j \mathbf{X}_{j-1}^\top M_n^{-1} \boldsymbol{\delta}(x), \qquad x \in \mathbb{R}.$$

We point out that for any square-integrable $a$,

$$a * \mathbb{J}_n(x) = \frac{1}{n - p_n} \sum_{j=p_n+1}^{n} \varepsilon_j \mathbf{X}_{j-1}^\top M_n^{-1} E[\mathbf{X}_0 a(x - Y_1)], \qquad x \in \mathbb{R}.$$

THEOREM 2. *Let (C), (I), (F), (S) and (1.8) hold and $p_n \to \infty$. Then, for each square-integrable $a$, the sequence $n^{1/2} a * \mathbb{J}_n$ is tight in $C_0(\mathbb{R})$ and $\|a * (\mathbb{H}_n - \mathbb{J}_n)\| = o_p(n^{-1/2})$.*

PROOF. Since $\mu_{n,i}(x) = E[X_0 k_{b_n}(x - Y_i)]$ equals $E[X_{1-i} k_{b_n}(x - Y_1)]$, we obtain that $\mathbb{H}_n = \hat{\Delta}^\top \mu_n$, where $\mu_n(x) = E[\mathbf{X}_0 k_{b_n}(x - Y_1)]$. Let us set

$$\tilde{\Delta} = M_n^{-1} \frac{1}{n - p_n} \sum_{j=p_n+1}^{n} \mathbf{X}_{j-1} \varepsilon_j.$$

By the results in Section 6, we have, with $v_n = k_{b_n} * \psi_1$ and $u_n = k_{b_n} * \psi_0$, that

$$\mu_{n,i}(x) = \mathbf{1}[i = \tau] E[v_n(x - Z_0)] + E[X_0 u_n(x - Z_i)].$$

Since $\|k_{b_n} * \psi_i - \psi_i\|_2 \to 0$ for $i = 0, 1$ and $\|k_{b_n} * \psi_0' - \psi_0'\|_2 \to 0$, we obtain from Lemma 6, applied with $a_n = k_{b_n}$, that

$$\sum_{i=1}^{\infty} \|\mu_{n,i} - \delta_i\|_2^2 \to 0 \quad \text{and} \quad \sum_{i=1}^{\infty} \|\delta_i\|_2^2 < \infty.$$

From this, we obtain that $\|\mu_n\|_2 = O(1)$. This shows that

(7.2) $\quad \|\mathbb{H}_n - \tilde{\Delta}^\top \mu_n\|_2 = \|(\hat{\Delta} - \tilde{\Delta})^\top \mu_n\|_2 \leq |\hat{\Delta} - \tilde{\Delta}| \|\mu_n\|_2 = o_p(n^{-1/2}).$

A martingale argument and straightforward calculations show that

$$\begin{aligned}(n - p_n) E[\mathbb{J}_n^2(x)] &= E[\varepsilon_0^2] E[(\mathbf{X}_0^\top M_n^{-1} \boldsymbol{\delta}(x))^2] \\ &= E[\varepsilon_0^2] E[\boldsymbol{\delta}(x)^\top M_n^{-1} \mathbf{X}_0 \mathbf{X}_0^\top M_n^{-1} \boldsymbol{\delta}(x)] \\ &= E[\varepsilon_0^2] \boldsymbol{\delta}(x)^\top M_n^{-1} M_n M_n^{-1} \boldsymbol{\delta}(x).\end{aligned}$$

This shows that

$$(n - p_n) E[\mathbb{J}_n^2(x)] \leq E[\varepsilon_0^2] K \sum_{i=1}^\infty \delta_i^2(x).$$

Since $\sum_{i=1}^\infty \delta_i^2$ is integrable, $n^{1/2} a * \mathbb{J}_n$ is tight by the results in Section 4. Since $\mu_{n,i} = k_{b_n} * \delta_i$, we find that $a * (\tilde{\Delta}^\top \mu_n) = k_{b_n} * a * \mathbb{J}_n$. Thus, by the tightness of $n^{1/2} a * \mathbb{J}_n$, we obtain that $\|a * (\tilde{\Delta}^\top \mu_n) - a * \mathbb{J}_n\| = o_p(n^{-1/2})$. This and (7.2) establish $n^{1/2} \|a * (\mathbb{H}_n - \mathbb{J}_n)\| = o_p(1)$. $\square$

Now, let us look at the case of parametric autocorrelation coefficients as described in Section 2. We then have $\varrho_i = r_i(\vartheta)$ and $\hat{\varrho}_i = r_i(\hat{\vartheta})$. We assume that (R1) and (R2) hold. This gives the expansion

$$R_n = \sum_{i=1}^{p_n} (r_i(\hat{\vartheta}) - r_i(\vartheta) - (\hat{\vartheta} - \vartheta)^\top \dot{r}_i(\vartheta))^2 = o_p(n^{-1}).$$

Fix a square-integrable $a$. Under (C), (S) and (F), we have

$$\sum_{i=1}^\infty \|a * \mu_{n,i} - a * \delta_i\|^2 \leq \|a\|_2^2 \sum_{i=1}^\infty \|\mu_{n,i} - \delta_i\|_2^2 \to 0$$

and

$$\sum_{i=1}^\infty \|a * \delta_i\|^2 \leq \|a\|_2^2 \sum_{i=1}^\infty \|\delta_i\|_2^2 < \infty.$$

Using the Cauchy–Schwarz inequality, we find that

$$\left\| \sum_{i=1}^{p_n} (r_i(\hat{\vartheta}) - r_i(\vartheta) - (\hat{\vartheta} - \vartheta)^\top \dot{r}_i(\vartheta)) a * \mu_{n,i} \right\|^2 \leq R_n \sum_{i=1}^\infty \|a * \mu_{n,i}\|^2 = o_p(n^{-1})$$

and

$$\left\| \sum_{i=1}^{p_n} \dot{r}_i(\vartheta) a * \mu_{n,i} - \sum_{i=1}^\infty \dot{r}_i(\vartheta) a * \delta_i \right\|^2$$
$$\leq \sum_{i=1}^\infty |\dot{r}_i(\vartheta)|^2 \left( \sum_{i=1}^{p_n} \|a * \mu_{n,i} - a * \delta_i\|^2 + \sum_{i=p_n+1}^\infty \|a * \delta_i\|^2 \right) \to 0,$$

provided $p_n \to \infty$. This shows that under (C), (I), (F), (R1), (R2) and (S), we have

$$\left\| a * \mathbb{H}_n - (\hat{\vartheta} - \vartheta)^\top \sum_{i=1}^\infty \dot{r}_i(\vartheta) a * \delta_i \right\| = o_p(n^{-1/2}).$$

Since $a * \delta_i(x) = E[X_0 a(x - Y_i)]$, we have the following result.

THEOREM 3. *Suppose that* (C), (I), (F), (R1), (R2) *and* (S) *hold and that* $\hat{\varrho}_i = r_i(\hat{\vartheta})$ *and* $\varrho_i = r_i(\vartheta)$. *Let* $p_n \to \infty$. *Then* $\|a * \mathbb{H}_n - (\hat{\vartheta} - \vartheta)^\top A\| = o_p(n^{-1/2})$, *where*

$$A(x) = \sum_{i=1}^\infty \dot{r}_i(\vartheta) E[X_0 a(x - Y_i)], \qquad x \in \mathbb{R}.$$

If $\dot{r}_i(\vartheta) = 0$ for all $i > p$, as is the case in the AR($p$) model, then the requirement that $p_n \to \infty$ can be relaxed to $p_n = p$.

**8. Behavior of the residuals.** In this section, we study how close the residuals are to the actual innovations. Recall that $\hat{\Delta} = (\hat{\varrho}_1 - \varrho_1, \ldots, \hat{\varrho}_{p_n} - \varrho_{p_n})^\top$ and $\mathbf{X}_{j-1} = (X_{j-1}, \ldots, X_{j-p_n})^\top$. Note that condition (R) is equivalent to $|\hat{\Delta}|^2 = O_p(q_n n^{-1})$. Under (I), we also have

$$\overline{\mathbf{X}} = \frac{1}{n - p_n} \sum_{j=p_n+1}^n \mathbf{X}_{j-1} = O_p(p_n^{1/2} n^{-1/2}).$$

This follows since we have

$$(8.1) \qquad (n - p_n) E\left[\left(\frac{1}{n - p_n} \sum_{j=p_n+1}^n X_{j-i}\right)^2\right] \le C E[X_0^2]$$

for some constant $C$ independent of $n$ and $i$. Thus, we derive

$$(8.2) \qquad \hat{\Delta}^\top \overline{\mathbf{X}} = O_p(p_n^{1/2} q_n^{1/2} n^{-1}).$$

The residuals can be expressed as

$$\hat{\varepsilon}_j = X_j - \sum_{i=1}^{p_n} \hat{\varrho}_i X_{j-i} = \varepsilon_j - \sum_{i=1}^{p_n} (\hat{\varrho}_i - \varrho_i) X_{j-i} + \sum_{i > p_n} \varrho_i X_{j-i} = \hat{\varepsilon}_j^* + \sum_{i > p_n} \varrho_i X_{j-i},$$

where

$$(8.3) \qquad \hat{\varepsilon}_j^* = \varepsilon_j - \sum_{i=1}^{p_n} (\hat{\varrho}_i - \varrho_i) X_{j-i} = \varepsilon_j - \hat{\Delta}^\top \mathbf{X}_{j-1}.$$

LEMMA 7. *Suppose that* (I), (Q) *and* (R) *hold. Then*

$$\sum_{j=p_n+1}^{n} (\hat{\varepsilon}_j - \hat{\varepsilon}_j^*)^2 = O_p(n^{-2\zeta}), \tag{8.4}$$

$$\sum_{j=p_n+1}^{n} (\hat{\varepsilon}_j^* - \varepsilon_j)^2 = O_p(p_n q_n), \tag{8.5}$$

$$\frac{1}{n-p_n} \sum_{j=p_n+1}^{n} (\hat{\varepsilon}_j - \varepsilon_j) = O_p(n^{-1/2-\zeta}) + O_p(p_n^{1/2} q_n^{1/2} n^{-1}). \tag{8.6}$$

*If the innovations have a finite moment of order $\xi \geq 2$, then*

$$\max_{p_n < j \leq n} |\hat{\varepsilon}_j - \varepsilon_j| = O_p(n^{-\zeta}) + o_p(p_n^{1/2} q_n^{1/2} n^{-1/2+1/\xi}). \tag{8.7}$$

PROOF. It follows from the Cauchy–Schwarz inequality that

$$(\hat{\varepsilon}_j^* - \varepsilon_j)^2 \leq \sum_{i=1}^{p_n} (\hat{\varrho}_i - \varrho_i)^2 \sum_{i=1}^{p_n} X_{j-i}^2. \tag{8.8}$$

From this bound, assumption (R) and the fact that $E[X_0^2] < \infty$, we obtain

$$\sum_{j=p_n+1}^{n} (\hat{\varepsilon}_j^* - \varepsilon_j)^2 = O_p(q_n n^{-1}) O_p(p_n n) = O_p(p_n q_n). \tag{8.9}$$

It follows from the Minkowski inequality that the $L_2(P)$-norm of $\hat{\varepsilon}_j - \hat{\varepsilon}_j^* = \sum_{s>p_n} \varrho_s X_{j-s}$ is bounded by the $L_2(P)$-norm of $X_0$ times $\sum_{s>p_n} |\varrho_s|$. Thus,

$$E\left[\sum_{j=p_n+1}^{n} (\hat{\varepsilon}_j - \hat{\varepsilon}_j^*)^2\right] \leq n E[X_0^2] \left(\sum_{s>p_n} |\varrho_s|\right)^2 = O(n^{-2\zeta}),$$

which implies (8.4). It follows from (8.4) that

$$\max_{p_n < j \leq n} |\hat{\varepsilon}_j - \hat{\varepsilon}_j^*| = O_p(n^{-\zeta}), \tag{8.10}$$

$$\frac{1}{n-p_n} \sum_{j=p_n+1}^{n} (\hat{\varepsilon}_j - \hat{\varepsilon}_j^*) = O_p(n^{-1/2-\zeta}). \tag{8.11}$$

Indeed, the square of the left-hand side of (8.10) is bounded by $R_n$, the left-hand side of (8.4), while the squared error term of (8.11) is bounded by $R_n/(n-p_n)$. Thus, (8.6) follows since, by (8.2), we have

$$\frac{1}{n-p_n} \sum_{j=p_n+1}^{n} (\hat{\varepsilon}_j^* - \varepsilon_j) = -\hat{\Delta}^\top \overline{\mathbf{X}} = O_p(p_n^{1/2} q_n^{1/2} n^{-1}). \tag{8.12}$$

The additional moment assumption on the innovations gives $E[|X_0|^\xi] < \infty$. From this, we obtain that $\max_{1 \leq j \leq n} |X_j| = o_p(n^{1/\xi})$. Indeed, for each $\eta > 0$,

$$P\left(\max_{1 \leq j \leq n} |X_j| > \eta n^{1/\xi}\right) \leq \sum_{j=1}^{n} P(|X_j| > \eta n^{1/\xi}) \leq \eta^{-\xi} E[X_0^\xi \mathbf{1}[|X_0| > \eta n^{1/\xi}]].$$

It follows from this, inequality (8.8) and assumption (R) that

$$(8.13) \quad \max_{p_n < j \leq n} |\hat{\varepsilon}_j^* - \varepsilon_j|^2 \leq p_n \sum_{i=1}^{p_n} (\hat{\varrho}_i - \varrho_i)^2 \max_{1 \leq j \leq n} |X_j|^2 = o_p(p_n q_n n^{-1+2/\xi}).$$

Combining (8.10) and (8.13), we obtain (8.7). □

LEMMA 8. *Suppose that* (I), (Q) *and* (R) *hold. Let $a_n$ be a sequence of functions with bounded integrable derivatives up to order two such that $\|a'_n\| = O(1)$ and $\|a''_n\| = o(p_n^{-1} q_n^{-1} n^{1/2})$. Then*

$$(8.14) \quad \sup_{x \in \mathbb{R}} \left| \frac{1}{n - p_n} \sum_{j=p_n+1}^{n} (a_n(x - \hat{Y}_j) - a_n(x - Y_j) + \hat{\Delta}^\top \mathbf{X}_{j-1} a'_n(x - Y_j)) \right|$$
$$= o_p(n^{-1/2}).$$

*If, further, $p_n q_n / n \to 0$ and $\|a''_n\|_2 = o(p_n^{-1/2} q_n^{-1/2} n^{1/2})$, then*

$$(8.15) \quad \sup_{x \in \mathbb{R}} \left| \frac{1}{n - p_n} \sum_{j=p_n+1}^{n} (a_n(x - \hat{\varepsilon}_j) - a_n(x - \varepsilon_j)) \right| = o_p(n^{-1/2}).$$

PROOF. Note that (8.4) implies

$$(8.16) \quad Q_n = \frac{1}{n - p_n} \sum_{j=p_n+1}^{n} |\hat{\varepsilon}_j - \hat{\varepsilon}_j^*| = O_p(n^{-\zeta - 1/2}),$$

while (8.3) and (8.5) imply

$$T_n = \frac{1}{n - p_n} \sum_{j=p_n+1}^{n} (\hat{\varepsilon}_j^* - \varepsilon_j)^2 = \frac{1}{n - p_n} \sum_{j=p_n+1}^{n} |\hat{\Delta}^\top \mathbf{X}_{j-1}|^2$$
$$(8.17) \quad = O_p(p_n q_n n^{-1}).$$

The expression following the supremum in (8.14) can be written as $|r_n(x)|$, where

$$r_n(x) = \frac{1}{n - p_n} \sum_{j=p_n+1}^{n} (a_n(x - \hat{Y}_j) - a_n(x - Y_j) + \hat{\Delta}^\top \mathbf{X}_{j-1} a'_n(x - Y_j)).$$

Define $r_n^*$ as $r_n$, but with $\hat{Y}_j = X_j - \hat{\varepsilon}_j$ replaced by $X_j - \hat{\varepsilon}_j^*$. Then

$$\|r_n - r_n^*\| \leq \|a_n'\| Q_n = O_p(n^{-\zeta-1/2}\|a_n'\|).$$

A Taylor expansion yields the bound

$$\|r_n^*\| \leq \|a_n''\| T_n = O_p(p_n q_n n^{-1}\|a_n''\|).$$

This establishes (8.14). The same arguments yield

$$\sup_{x\in\mathbb{R}}\left|\frac{1}{n-p_n}\sum_{j=p_n+1}^{n}(a_n(x-\hat{\varepsilon}_j) - a_n(x-\varepsilon_j) - \hat{\Delta}^\top \mathbf{X}_{j-1} a_n'(x-\varepsilon_j))\right| = o_p(n^{-1/2}).$$

In view of (8.2), we have

$$\|\hat{\Delta}^\top \overline{\mathbf{X}} a_n' * f\| \leq |\hat{\Delta}^\top \overline{\mathbf{X}}|\|a_n' * f\| = O_p(p_n^{1/2} q_n^{1/2} n^{-1}\|a_n'\|) = o_p(n^{-1/2}).$$

Result (8.15) now follows if we can show that $\|\hat{\alpha}_n\| = o_p(q_n^{-1/2})$ for

$$\hat{\alpha}_n(x) = \frac{1}{n-p_n}\sum_{j=p_n+1}^{n}\mathbf{X}_{j-1}(a_n'(x-\varepsilon_j) - E[a_n'(x-\varepsilon_j)]), \qquad x\in\mathbb{R}.$$

It follows from Fubini's theorem that $\hat{\alpha}_n = a_n'' * W_n$ with

$$W_n(x) = \frac{1}{n-p_n}\sum_{j=p_n+1}^{n}\mathbf{X}_{j-1}(\mathbf{1}[\varepsilon_j \leq x] - F(x)).$$

Thus, $\|\hat{\alpha}_n\| \leq \|a_n''\|_2 \|W_n\|_2$. Since

$$(n-p_n)E[\|W_n\|_2^2] = E[|\mathbf{X}_0|^2]\int F(x)(1-F(x))\,dx = O(p_n),$$

we obtain $\|\hat{\alpha}_n\| = O_p(p_n^{1/2} n^{-1/2}\|a_n''\|_2) = o_p(q_n^{-1/2})$.  $\square$

**9. Estimating the innovation density $f$.** The kernel estimator based on the residuals is

$$\hat{f}(x) = \frac{1}{n-p_n}\sum_{j=p_n+1}^{n} k_{b_n}(x-\hat{\varepsilon}_j), \qquad x\in\mathbb{R}.$$

In this section, we study convergence of $\hat{f}$ in the space $L_2$ and of functionals of the form $a * \hat{f}$ in the space $C_0(\mathbb{R})$.

Let $\tilde{f}$ denote the kernel estimator based on the actual innovations $\varepsilon_{p_n+1},\ldots,\varepsilon_n$,

$$\tilde{f}(x) = \frac{1}{n-p_n}\sum_{j=p_n+1}^{n} k_{b_n}(x-\varepsilon_j), \qquad x\in\mathbb{R}.$$

The first result is known.

LEMMA 9. *Suppose that the kernel $k$ is square-integrable and of type $(m,2)$. Let $f$ be $L_2$-smooth of order $r \leq m$. Then*
$$\|\tilde{f} - f\|_2 = O_p(b_n^{-1/2} n^{-1/2}) + o(b_n^r).$$

PROOF. It is well known that $E[\tilde{f}(x)] = f * k_{b_n}(x)$ and
$$(n - p_n) E[\|\tilde{f} - f * k_{b_n}\|_2^2] \leq \|k_{b_n}^2 * f\|_1 \leq b_n^{-1} \|k^2\|_1.$$
Thus, $\|\tilde{f} - f * k_{b_n}\|_2 = O_p(b_n^{-1/2} n^{-1/2})$. By Corollary 2, $\|f * k_{b_n} - f\|_2 = o(b_n^r)$. □

LEMMA 10. *Suppose that* (I), (Q), (R), (F) *and* (K) *hold. Then*
$$\|\hat{f} - \tilde{f}\|_2 = O_p(p_n q_n b_n^{-5/2} n^{-1}) + O_p(n^{-\zeta - 1/2} b_n^{-3/2}).$$

PROOF. Let $\hat{\varepsilon}_j^*$ be as in (8.3). Let $\hat{f}^*$ denote the kernel estimator based on $\hat{\varepsilon}_{p_n+1}^*, \ldots, \hat{\varepsilon}_n^*$. With $Q_n$ as in (8.16), we find that
$$\|\hat{f} - \hat{f}^*\|_2^2 \leq \|\hat{f} - \hat{f}^*\|_1 \|\hat{f} - \hat{f}^*\| \leq \|k_{b_n}'\|_1 \|k_{b_n}'\| Q_n^2$$
and obtain, in view of (8.16), the rate
$$\|\hat{f} - \hat{f}^*\|_2 = O_p(b_n^{-3/2} n^{-\zeta - 1/2}).$$
The identity $\hat{\varepsilon}_j^* = \varepsilon_j - \hat{\Delta}^\top \mathbf{X}_{j-1}$ and a Taylor expansion yield $\hat{f}^* - \tilde{f} = \hat{\Delta}^\top \gamma_n + r_n$ with
$$\gamma_n(x) = \frac{1}{n - p_n} \sum_{j=p_n+1}^{n} \mathbf{X}_{j-1} k_{b_n}'(x - \varepsilon_j),$$
$$r_n(x) = \frac{1}{n - p_n} \sum_{j=p_n+1}^{n} \int_0^1 \int_0^1 (\hat{\Delta}^\top \mathbf{X}_{j-1})^2 t k_{b_n}''(x - \varepsilon_j + st \hat{\Delta}^\top \mathbf{X}_{j-1}) \, ds \, dt.$$
With $T_n$ as in (8.17), we obtain $\|r_n\| \leq \|k_{b_n}''\| T_n = O_p(p_n q_n b_n^{-3} n^{-1})$ and $\|r_n\|_1 \leq \|k_{b_n}''\|_1 T_n = O_p(p_n q_n b_n^{-2} n^{-1})$, and, consequently,
$$\|r_n\|_2^2 \leq \|r_n\| \|r_n\|_1 = O_p(p_n^2 q_n^2 b_n^{-5} n^{-2}).$$
Let $\bar{\gamma}_n = \overline{\mathbf{X}} k_{b_n}' * f$. Since $\|k_{b_n}' * f\|_2 = \|f' * k_{b_n}\|_2 \leq \|f'\|_2 \|k\|_1$, we obtain from (8.2) that
$$\|\hat{\Delta}^\top \bar{\gamma}_n\|_2 \leq |\hat{\Delta}^\top \overline{\mathbf{X}}| \|k_{b_n}' * f\|_2 = O_p(p_n^{1/2} q_n^{1/2} n^{-1}).$$
A martingale argument yields
$$(n - p_n) E[\|\gamma_n - \bar{\gamma}_n\|_2^2] \leq p_n E[|X_0^2|] \|(k_{b_n}')^2 * f\|_1 = O(p_n b_n^{-3}).$$
Thus, $\|\hat{\Delta}^\top (\gamma_n - \bar{\gamma}_n)\|_2 = O_p(p_n^{1/2} q_n^{1/2} b_n^{-3/2} n^{-1})$. The above imply the desired rate. □

THEOREM 4. *Suppose that* (I), (Q), (R), (F) *and* (K) *hold. Let* $a \in \mathcal{A}$ *and let* $a * f$ *be smooth of order* $r \leq m$. *Let the bandwidth satisfy* $nb_n^{2r} = O(1)$ *and* $p_n q_n b_n^{-1} n^{-1/2} \to 0$. *Then*

$$\|a * (\hat{f} - f) - \mathbb{A}_n\| = o_p(n^{-1/2}),$$

*where*

$$\mathbb{A}_n(x) = \frac{1}{n - p_n} \sum_{j=p_n+1}^{n} (a(x - \varepsilon_j) - E[a(x - \varepsilon_j)]), \qquad x \in \mathbb{R}.$$

PROOF. Let $\bar{f} = E[\tilde{f}] = f * k_{b_n}$. Since $a * f$ is smooth of order $r \leq m$ and $k$ is of type $(m, 1)$, Corollary 3 yields

$$\|a * \bar{f} - a * f\| = \|(a * f) * k_{b_n} - a * f\| = o(b_n^r) = o(n^{-1/2}).$$

We can write $a * (\tilde{f} - \bar{f}) = \mathbb{A}_n * k_{b_n}$. Since $n^{1/2} \mathbb{A}_n$ is tight in $C_0(\mathbb{R})$ by result (a) in Section 4, we obtain that $\|n^{1/2}(\mathbb{A}_n * k_{b_n} - \mathbb{A}_n)\| = o_p(1)$. In other words,

$$\|a * (\tilde{f} - \bar{f}) - \mathbb{A}_n\| = o_p(n^{-1/2}).$$

We can now calculate that

$$a * (\hat{f} - \tilde{f})(x) = \frac{1}{n - p_n} \sum_{j=p_n+1}^{n} (a_n(x - \hat{\varepsilon}_j) - a_n(x - \varepsilon_j)), \qquad x \in \mathbb{R},$$

with $a_n = a * k_{b_n}$. $a_n$ is then twice differentiable with $a_n' = a' * k_{b_n}$ and $a_n'' = a' * k_{b_n}'$. We have $\|a_n'\| \leq \|a'\| \|k_{b_n}\|_1 = O(1)$, $\|a_n''\| \leq \|a'\| \|k_{b_n}'\|_1 = O(b_n^{-1})$ and $\|a_n''\|_2^2 \leq \|a_n''\| \|a_n''\|_1 \leq \|a_n''\| \|k_{b_n}'\|_1 \|a'\|_1 = O(b_n^{-2})$. In view of $p_n q_n b_n^{-1} n^{-1/2} \to 0$, Lemma 8 yields

$$\|a * (\hat{f} - \tilde{f})\| = o_p(n^{-1/2}).$$

The desired result follows from the above. □

**10. Estimating the density $g$.** The kernel estimator based on the estimated versions $\hat{Y}_j = X_j - \hat{\varepsilon}_j$ of the $Y_j = X_j - \varepsilon_j$ is

$$\hat{g}(x) = \frac{1}{n - p_n} \sum_{j=p_n+1}^{n} k_{b_n}(x - \hat{Y}_j), \qquad x \in \mathbb{R}.$$

In this section, we study convergence of $\hat{g}$ in the space $L_2$ and of functionals of the form $a * \hat{g}$ in the space $C_0(\mathbb{R})$. Let $\tilde{g}$ denote the kernel estimator based on $Y_{p_n+1}, \ldots, Y_n$,

$$\tilde{g}(x) = \frac{1}{n - p_n} \sum_{j=p_n+1}^{n} k_{b_n}(x - Y_j), \qquad x \in \mathbb{R}.$$

We first give an analogue of Lemma 9.

LEMMA 11. *Suppose that* (C) *and* (S) *hold. Let the kernel $k$ be square-integrable and of type $(m,2)$. Let $f$ belong to $\mathcal{A}_1 \cap \mathcal{A}_2$ and have finite mean. Let $g$ be $L_2$-smooth of order $r$ with $r \leq m$. Then*

$$\|\tilde{g} - g\|_2 = O_p(b_n^{-1/2} n^{-1/2}) + o(b_n^r).$$

PROOF. By Corollary 2, we have $\|g * k_{b_n} - g\|_2 = o(b_n^r)$. It remains to show that

(10.1) $$\|\tilde{g} - g * k_{b_n}\|_2 = O_p(b_n^{-1/2} n^{-1/2}).$$

Recall the notation $\tau = \inf\{s \geq 1 : \varphi_s \neq 0\}$. We can write $Y_j = \varphi_\tau \varepsilon_{j-\tau} + Z_j$ with $Z_j = \sum_{s > \tau} \varphi_s \varepsilon_{j-s}$. Let $a_n = k_{b_n} * \psi_0$, where $\psi_0$ is the density of $\varphi_\tau \varepsilon_0$. We can then express $\tilde{g} - g * k_{b_n}$ as the sum $T_1 + k_{b_n} * T_2$ with

$$T_1(x) = \frac{1}{n - p_n} \sum_{j=p_n+1}^{n} (k_{b_n}(x - Y_j) - a_n(x - Z_j)),$$

$$T_2(x) = \frac{1}{n - p_n} \sum_{j=p_n+1}^{n} (\psi_0(z - Z_j) - E[\psi_0(x - Z_j)]).$$

Using a martingale argument, we obtain $(n - p_n) E[\|T_1\|_2^2] \leq \|k_{b_n}^2 * g\|_1 = O(b_n^{-1})$ and thus $\|T_1\|_2 = O_p(b_n^{-1/2} n^{-1/2})$. Since $f$ belongs to $\mathcal{A}_1 \cap \mathcal{A}_2$, so does $\psi_0$. Thus, $n^{1/2} T_2$ is tight by result (b) in Section 4, applied with $A = \psi_0$ and $\xi_j = Z_j$. This shows that $\|T_2 * k_{b_n}\|_2^2 \leq \|T_2\|_2^2 \|k_{b_n}\|_1 \leq \|T_2\| \|T_2\|_1 \|k\| = O_p(n^{-1/2})$. This finishes the proof of (10.1). $\square$

Let us define functions $\mu_n$ and $\mu'_n$ by

$$\mu_n(x) = E[\mathbf{X}_0 k_{b_n}(x - Y_1)] \quad \text{and} \quad \mu'_n(x) = E[\mathbf{X}_0 k'_{b_n}(x - Y_1)].$$

We now give analogues of Lemma 10 and Theorem 4.

LEMMA 12. *Suppose that* (C), (I), (Q), (R), (S), (F) *and* (K) *hold. Then*

$$\|\hat{g} - \tilde{g} + \hat{\Delta}^\top \mu'_n\|_2 = O_p(p_n q_n b_n^{-5/2} n^{-1}) + O_p(n^{-\zeta - 1/2} b_n^{-3/2}).$$

PROOF. Let $\hat{g}^*$ denote the kernel estimator based on $\hat{Y}^*_{p_n+1}, \ldots, \hat{Y}^*_n$ with

$$\hat{Y}^*_j = X_j - \hat{\varepsilon}^*_j = Y_j + \hat{\Delta}^\top \mathbf{X}_{j-1}.$$

As in the proof of Lemma 10, we find that

$$\|\hat{g} - \hat{g}^*\|_2 = O_p(n^{-\zeta - 1/2} b_n^{-3/2}) \quad \text{and}$$

$$\|\hat{g}^* - \tilde{g} + \hat{\Delta}^\top \hat{\mu}'_n\|_2 = O_p(p_n q_n b_n^{-5/2} n^{-1}),$$

where

$$\hat{\mu}'_n(x) = \frac{1}{n - p_n} \sum_{j=p_n+1}^{n} \mathbf{X}_{j-1} k'_{b_n}(x - Y_j), \qquad x \in \mathbb{R}.$$

Note that $\|k'_{b_n}\| = O(b_n^{-2})$ and $\|k'_{b_n}\| = O(b_n^{-1})$. Thus, it follows from Lemma 5, applied with $a_n = k'_{b_n}$, that

$$\int E[\|\hat{\mu}'_n(x) - E[\hat{\mu}'_n(x)]\|^2]\, dx = O(p_n b_n^{-3} n^{-1}).$$

Since $\mu'_n(x) = E[\hat{\mu}'_n(x)]$, we see that

$$\|\hat{\Delta}^\top (\hat{\mu}'_n - \mu'_n)\|_2 = O_p(p_n^{1/2} q_n^{1/2} b_n^{-3/2} n^{-1}).$$

The above rates yield the desired result. $\square$

THEOREM 5. *Suppose that* (C), (I), (Q), (R), (S), (F) *and* (K) *hold. Let $a \in \mathcal{A}$ and let $a * g$ be smooth of order $r$ with $r \leq m$. Let the bandwidth satisfy $n b_n^{2r} = O(1)$ and $p_n q_n b_n^{-1} n^{-1/2} \to 0$. Then*

$$\|a * (\hat{g} - g) - \mathbb{K}_n + a' * (\hat{\Delta}^\top \mu_n)\| = o_p(n^{-1/2}),$$

where

$$\mathbb{K}_n(x) = \frac{1}{n - p_n} \sum_{j=p_n+1}^{n} (a(x - Y_j) - E[a(x - Y_j)]), \qquad x \in \mathbb{R}.$$

PROOF. Set $\bar{g} = E[\tilde{g}] = g * k_{b_n}$. Since $a * g$ is smooth of order $r$ and the kernel $k$ is of type $(m, 1)$ with $m \geq r$, we obtain from Corollary 3 that

$$\|a * \bar{g} - a * g\| = \|(a * g) * k_{b_n} - a * g\| = o(b_n^r) = o(n^{-1/2}).$$

Simple calculations yield $a * (\tilde{g} - \bar{g}) = \mathbb{K}_n * k_{b_n}$. Since $a$ belongs to $\mathcal{A}_1 \cap \mathcal{A}_2$ and $f$ has finite mean, it follows from (S) and result (b) in Section 4 that $n^{1/2} \mathbb{K}_n$ is tight in $C_0(\mathbb{R})$. Consequently, $\|n^{1/2}(\mathbb{K}_n * k_{b_n} - \mathbb{K}_n)\| = o_p(1)$. In other words,

$$\|a * (\tilde{g} - \bar{g}) - \mathbb{K}_n\| = o_p(n^{-1/2}).$$

With $a_n = a * k_{b_n}$, one verifies that

$$a * (\hat{g} - \tilde{g})(x) = \frac{1}{n - p_n} \sum_{j=p_n+1}^{n} (a_n(x - \hat{Y}_j) - a_n(x - Y_j)), \qquad x \in \mathbb{R}.$$

Now, let

$$\hat{\mu}_n(x) = \frac{1}{n - p_n} \sum_{j=p_n+1}^{n} \mathbf{X}_{j-1} k_{b_n}(x - Y_j), \qquad x \in \mathbb{R}.$$

Since $\|a'_n\| = O(1)$, $\|a''_n\| = O(b_n^{-1})$ and $\|a''_n\|_2 = O(b_n^{-1})$, as shown in the proof of Theorem 4, and since $p_n q_n b_n^{-1} n^{-1/2} \to 0$, we obtain from Lemma 8 and $a'_n = a' * k_{b_n}$ that

$$\|a * (\hat{g} - \tilde{g}) + a' * (\hat{\Delta}^\top \hat{\mu}_n)\| = o_p(n^{-1/2}).$$

It follows from Lemma 5, $\|k_{b_n}\| = O(b_n^{-1})$ and $\|k_{b_n}\|_1 = O(1)$ that

$$\int E[\|\hat{\mu}_n(x) - E[\hat{\mu}_n(x)]\|^2]\,dx = O_p(p_n b_n^{-1} n^{-1}).$$

Since $\mu_n(x) = E[\hat{\mu}_n(x)]$, we find that

$$\|a' * \hat{\Delta}^\top (\hat{\mu}_n - \mu_n)\| \le \|a'\|_2 |\hat{\Delta}| \|\hat{\mu}_n - \mu_n\|_2 = O_p(p_n^{1/2} q_n^{1/2} b_n^{-1/2} n^{-1})$$
$$= o_p(n^{-1/2}).$$

The desired result follows from the above. $\square$

**Acknowledgments.** We thank an Associate Editor and two referees for suggestions that led us to completely rewrite this paper. We had originally introduced a more complicated $n^{1/2}$-consistent density estimator that involved an increasing number of convolutions; one referee suggested that a single convolution should suffice.

DEPARTMENT OF MATHEMATICAL SCIENCES
BINGHAMTON UNIVERSITY
BINGHAMTON, NEW YORK 13902-6000
USA
E-MAIL: anton@math.binghamton.edu

MATHEMATICAL INSTITUTE
UNIVERSITY OF COLOGNE
50931 COLOGNE
GERMANY
E-MAIL: wefelm@math.uni-koeln.de